\documentclass[12pt]{article}
\usepackage[colorlinks,
            linkcolor=red,
            anchorcolor=blue,
            citecolor=blue
            ]{hyperref}
\usepackage{mathrsfs,amssymb,amsfonts}
\usepackage{amsmath,amssymb,latexsym,color}
\usepackage[mathscr]{eucal}
\usepackage[numbers,sort&compress]{natbib}
\usepackage{graphicx}
\usepackage{lipsum}
\def\thefootnote{\fnsymbol{footnote}}
\newtheorem{thm}{Theorem}[section]

\newtheorem{prop}[thm]{Proposition}
\newtheorem{lemma}[thm]{Lemma}
\newtheorem{cor}[thm]{Corollary}
\newtheorem{example}[thm]{Example}

\newtheorem{problem}[thm]{Problem}
\newtheorem{remark}[thm]{Remark}
\newtheorem{Notation}[thm]{Notation}

\newcommand{\proof}{{\it Proof.\quad}}
\newcommand{\qed}{\hfill\Box\medskip}
\renewcommand{\thefootnote}{\arabic{footnote}}

\newcommand{\sdim}{{\rm sdim}}
\usepackage{CJK}
\begin{document}
\begin{CJK*}{GBK}{song}
\renewcommand{\abovewithdelims}[2]{
\genfrac{[}{]}{0pt}{}{#1}{#2}}
%%%%%%%%%%%%%%%%%%%%%%%%%%%%%%%%%%%%%%%%%%%%%%%%%%%%%%%%%%%%%%%%%%%%%%%%%%%%%%%%%%%%%%%%
%%%%%%%%%%%%%%%%%%%%%%%%%%%%%%%%%%%%%%%%%%%%%%%%%%%%%%%%%%%%%%%%%%%%%%%%%%%%%%%%%%%%%%%%

\title{\bf Strong metric dimensions for power graphs of  finite groups}

\author{{Xuanlong Ma$^{1,}$\footnote{Corresponding author} ~and Liangliang Zhai$^{2}$}
\\
{\small\em $^1$School of Science, Xi'an Shiyou University, Xi'an  710065, China}\\
{\small\em $^2$School of Mathematics, Northwest University,
Xi'an 710127, China}\\
}

 \date{}
 \maketitle
\newcommand\blfootnote[1]{%
\begingroup
\renewcommand\thefootnote{}\footnote{#1}%
\addtocounter{footnote}{-1}%
\endgroup
}
\begin{abstract}
Let $G$ be a finite group.
The order supergraph of $G$
is the graph with vertex
set $G$, and two distinct vertices $x,y$ are adjacent if
$o(x)\mid o(y)$ or $o(y)\mid o(x)$. The enhanced power graph of $G$ is the graph whose vertex set is $G$, and two distinct vertices are adjacent if they generate a cyclic subgroup.
The reduced power graph of $G$ is the graph with vertex
set $G$, and two distinct vertices $x,y$  are adjacent if
$\langle x\rangle \subset \langle y\rangle$ or $\langle y\rangle \subset \langle x\rangle$. In this paper,
we characterize the strong metric dimension of
the order supergraph, the enhanced power graph and
the reduced power graph of a finite group.


\end{abstract}

\blfootnote{E-mail address: xuanlma@mail.bnu.edu.cn (X. Ma), zhailiang111@126.con (L. Zhai)}

\section{Introduction}
All graphs considered in this paper are finite, undirected, with no loops and no multiple edges. Let $\Gamma$ be a graph.
The vertex set of $\Gamma$ is denoted by $V(\Gamma)$.
Let $x,y,z\in V(\Gamma)$.
The {\em distance} between $x$ and $y$ in $\Gamma$, denoted by $d(x,y)$, is the length
of a shortest path from $x$ to $y$.
The {\em diameter} of $\Gamma$ is the greatest distance between any two vertices.
We say that $z$ {\em strongly resolves} $x$ and $y$ if there exists a shortest path from $z$ to $x$ containing $y$, or a  shortest path from $z$ to $y$ containing $x$. A subset $S$ of $V(\Gamma)$ is a {\em strong
resolving set} of $\Gamma$ if every pair of vertices of $\Gamma$ is strongly resolved by some vertex in $S$. The smallest cardinality of a strong resolving set of $\Gamma$,
denoted by $\sdim(\Gamma)$,
is called the {\em strong metric dimension} of $\Gamma$.

In the 1970s, the metric dimension of a graph
was introduced independently by
Harary and Melter \cite{Ham} and Slater \cite{Sl}.
In 2004, Seb\H{o} and Tannier \cite{ST04} introduced
the strong metric dimension of a graph and presented some applications of strong resolving sets to combinatorial searching.
The problem of computing strong metric dimension is NP-hard \cite{KYR3}.
Some theoretical results, computational approaches and recent results on strong metric dimension can be found in \cite{KKCM14}.

Graphs associated with groups and other algebraic structures have been actively investigated, since they have valuable applications
(cf. \cite{Knew,KeR}) and are related to automata theory (cf. \cite{K,K1}).
The {\em undirected power graph} $\mathcal{P}(G)$ of a finite group $G$ has vertex set $G$ and two distinct elements are adjacent if one is a power of the other.
The concepts of power graph
and undirected power graph were first introduced by Kelarev and Quinn \cite{KQ00} and Chakrabarty {\em et al.} \cite{CGS}, respectively. The metric dimension and the strong metric dimension of a power
graph were studied in \cite{FMW} and \cite{Mas}, respectively.
In recent years, the study of power graphs has been growing, see, for example, \cite{Cam19,B2,CGh,MWW,MRS}.
Also, see \cite{AKC} for a survey of results and open problems on power graphs.

Let $G$ be a finite group.
The {\em enhanced power graph} $\mathcal{P}_E(G)$ of $G$ is the graph whose vertex set is $G$, and two distinct vertices are adjacent if they generate a cyclic subgroup of $G$.
In order to measure how close the power graph is to the commuting graph, Aalipour  {\em et al.} \cite{acam} introduced the enhanced power graph which lies in between.
Ma and She \cite{MS} characterized the metric dimension of an enhanced power graph.
See \cite{acam,Bh18,CM,CM2,PDK} for some more properties of the enhanced power graph.

The {\em order supergraph} $\mathcal{S}(G)$ of $\mathcal{P}(G)$ of
$G$ is a graph with vertex
set $G$, and two distinct vertices $x,y$ are adjacent if
$o(x)\mid o(y)$ or $o(y)\mid o(x)$,
where $o(x)$ and $o(y)$ are the orders of $x$ and $y$, respectively.
By the definition of an order supergraph,
we also call $\mathcal{S}(G)$ as the {\em order supergraph} of $G$.
In 2017, Hamzeh and Ashrafi \cite{HA17} called
this graph as the {\em main supergraph} of $G$ and studied its full automorphism group. Recently, Hamzeh and Ashrafi \cite{HA}
studied some properties of the order supergraph, and in particular, they showed that
$\mathcal{S}(G)=\mathcal{P}(G)$ if and only if $G$ is cyclic. Also, in \cite{HA19}, they investigated  Hamiltonianity,
Eulerianness and $2$-connectedness of this graph.

With an intention to avoid the complexity of edges in the power graphs, Rajkumar and Anitha \cite{RP1} introduced the {\em reduced power graph} $\mathcal{P}_{R}(G)$ of $G$, which is an undirected graph with vertex
set $G$, and two distinct vertices $x,y$  are adjacent if
$\langle x\rangle \subset \langle y\rangle$ or $\langle y\rangle \subset \langle x\rangle$. In other words, $\mathcal{P}_{R}(G)$  is the subgraph of $\mathcal{P}(G)$ obtained by deleting all edges $\{x,y\}$ with $\langle x\rangle=\langle y\rangle$, where $x$ and $y$ are two distinct elements of $G$. In \cite{RP1},
the authors studied the interplay between
the algebraic properties of a group and the graph theoretic properties of its reduced power graph. Recently, Anitha and Rajkumar \cite{AR19} characterized the groups with planar, toroidal and projective planar reduced power graphs.
Moreover, see \cite{AnR,RP2} for some more properties of this  graph.

According to the definitions as above,
for any finite group $G$,
$\mathcal{P}_{R}(G)$ is a spanning subgraph of $\mathcal{P}(G)$, and $\mathcal{P}(G)$ is a spanning subgraph of both $\mathcal{S}(G)$ and $\mathcal{P}_E(G)$. In this paper,
we characterize the strong metric dimension of
the order supergraph, the enhanced power graph and
the reduced power graph of a finite group.

\section{Preliminaries}\label{pre}
This section introduces some basic definitions and notations that are used throughout the paper.

Every group considered in this paper is finite.
We always use $e$ to denote the identity
element of the group under consideration.
Let $G$  be a group. The {\em order} of an element $x$ of $G$, denoted by $o(x)$,
is defined as
the cardinality of the cyclic subgroup $\langle x\rangle$. An element of order $2$ is called an {\em involution}.
The {\em exponent} of $G$, denoted by $\exp(G)$,
is defined as the least common multiple of the orders of all elements of $G$.
The set of orders of all elements of $G$ is denoted by $\pi_e(G)$.
A {\em maximal cyclic subgroup} of $G$ is a cyclic subgroup, which is not a proper subgroup of some cyclic
subgroup of $G$. The set of all maximal cyclic subgroups of $G$
is denoted by $\mathcal{M}_G$. Note
that $|\mathcal{M}_G|=1$ if and only if $G$ is cyclic. Denote by $\mathbb{Z}_n$ the cyclic group of order $n$.

A finite group is called a {\em $\mathcal{P}$-group} \cite{Dea} if every nontrivial element of the group has prime order. For example, the elementary abelian $p$-group $\mathbb{Z}_p^n$ is a $\mathcal{P}$-group where $p$ is a prime and $n\ge 1$, and the symmetric group $S_3$ on $3$ letters is also a $\mathcal{P}$-group.
A finite group is called a {\em CP-group} \cite{Hi} if
every nontrivial element of the group has prime power order. Clearly, both $p$-groups and $\mathcal{P}$-groups are also  CP-groups.

For $n\ge 2$, Johnson \cite[pp. 44--45]{Jon} defined the generalized quaternion group $Q_{4n}$ of order $4n$ by the presentation
\begin{equation}\label{q4n}
Q_{4n}=\langle x,y: x^n=y^2,x^{2n}=y^4=e,y^{-1}xy=x^{-1}\rangle.
\end{equation}
If $n=2$, then $Q_8$ is the usual quaternion group of order $8$.
Some basic properties of $Q_{4n}$ can be found in \cite{Gor}. We remark that $x^n$ is the unique involution of $Q_{4n}$. Also, it is easy to check that
\begin{equation}\label{q4n-0}
Q_{4n}=\langle x\rangle\cup \{x^iy: 1\le i \le 2n\}, ~o(x^iy)=4\text{ for each $1\le i \le 2n$}
\end{equation}
and
\begin{equation}\label{q4n-1}
\mathcal{M}_{Q_{4n}}=\{\langle x\rangle,\langle xy\rangle,
\ldots,\langle x^ny\rangle\},~~x^n\in \bigcap_{M\in \mathcal{M}_{4n}}M.
\end{equation}

Recall now the following elementary result.
\begin{thm}\label{a-p-group}{\rm (\cite[Theorem~5.4.10~(ii)]{Gor})}
\label{p-unique}
Let $p$ be a prime. Then a $p$-group having
a unique subgroup of order $p$ is either cyclic or generalized quaternion.
\end{thm}

Let $\Gamma$ be a graph and $x\in V(\Gamma)$.
The {\em closed neighborhood} of $x$ in $\Gamma$ is $$N_\Gamma[x]=\{y\in V(\Gamma): d(y,x)\le 1\}.$$
If the situation is unambiguous, we denote $N_\Gamma[x]$ simply by $N[x]$.
A subset of $V(\Gamma)$ is called
a {\em clique} if any two distinct vertices in this subset are adjacent
in $\Gamma$. The {\em clique number} of $\Gamma$, denoted by  $\omega(\Gamma)$, is the maximum cardinality of a clique in $\Gamma$.

For $x,y\in V(\Gamma)$, define a binary relation $x\approx y$  by the rule that $N[x]=N[y]$ in $\Gamma$.
Observe  that $\approx$ is an equivalence relation over $V(\Gamma)$.
Let $U(\Gamma)$ be a complete set of distinct representative elements for this equivalence relation.
The {\em reduced graph} of $\Gamma$, denoted by $\mathcal{R}_\Gamma$,
has the vertex set $U(\Gamma)$ and two vertices are adjacent if they are adjacent in $\Gamma$. Notice that
for two distinct equivalence classes $\mathcal{C}_1$ and $\mathcal{C}_2$, if there exist a vertex in $\mathcal{C}_1$ and a vertex in $\mathcal{C}_2$ which are adjacent in $\Gamma$, then each vertex in $\mathcal{C}_1$ and each vertex in $\mathcal{C}_2$ are adjacent in $\Gamma$.
As a result, $\mathcal{R}_\Gamma$ does not depend on the choice of representatives.

Ma {\em et al.} \cite{Mas} characterized the strong metric dimension of a graph with diameter two by the reduced graph of this graph.

\begin{thm}{\rm (\cite[Theorem 2.2]{Mas})}\label{Newthm0}
Let $\Gamma$ be a connected graph with order $n$ and diameter two. Then
$
\sdim(\Gamma)=n-\omega(\mathcal{R}_\Gamma).
$
\end{thm}

For a positive integer $n$, let
$n=p_1^{r_1}p_2^{r_2}\cdots p_m^{r_m}$
be its canonical factorization, that is,
$p_1,p_2,\ldots,p_m$ are pairwise distinct primes and $r_i\ge 1$ for $1\le i \le m$. Denote by $\Omega(n)$ the number of all prime factors of $n$ counted with multiplicity. Namely,
$$
\Omega(n)=\sum_{i=1}^mr_i.
$$

\section{Order supergraphs of power graphs}
This section characterizes the strong metric dimension of the order supergraph of a group.
Our main result is as follows.

\begin{thm}\label{mainth-or}
Let $G$ be a group of order $n$. Then
$$\sdim(\mathcal{S}(G))
=\left\{
                                  \begin{array}{ll}
                                  n-1, & \hbox{if $G$ is a $p$-group;} \\
                                  n-\Omega(n), & \hbox{if $G$ is a cyclic group and is not a $p$-group;} \\
                                  n-2, & \hbox{if $G$ is a CP-group and is not a $p$-group;} \\
                                  n-\lambda_G-1, & \hbox{otherwise,}
                                  \end{array}
                                \right.
$$
where $\lambda_G=\max\{\Omega(m): m\in \pi_e(G) \text{ and $m$ is not a prime power}\}$.
\end{thm}

Note that $\mathcal{S}(G)$ is complete if and only if $G$ is a $p$-group (see also \cite[Theorem 2.3]{HA}).
So, $\sdim(\mathcal{S}(G))=|G|-1$ if and only if $G$ is a $p$-group.
As a corollary of Theorem~\ref{mainth-or},
we can classify all groups $G$
whose order supergraphs have strong metric dimension $|G|-2$.

\begin{cor}
Let $G$ be a group of order $n$. Then
$\sdim(\mathcal{S}(G))=n-2$ if and only if
$G$ is isomorphic to either $\mathbb{Z}_{pq}$ or a CP-group with at least two distinct prime divisors, where $p,q$ are two distinct primes.
\end{cor}

By Theorem~\ref{mainth-or} and (\ref{q4n-0}), we determine the strong metric dimension of the order supergraph of a generalized quaternion group.

\begin{cor}
Let $Q_{4n}$ be the generalized quaternion group as presented in {\rm (\ref{q4n})}. Then
$$\sdim(\mathcal{S}(Q_{4n}))
=\left\{
                                  \begin{array}{ll}
                                    4n-1, & \hbox{if $n$ is a power of $2$;} \\
                                    4n-\Omega(2n)-1, &  \hbox{otherwise.}
                                  \end{array}
                                \right.
$$
\end{cor}

In the following, we aim to prove Theorem~\ref{mainth-or}.
For $x,y \in G$, denote by $\sim$ the equivalence relation defined by  $N[x]=N[y]$ in $\mathcal{S}(G)$.
As stated above, $\sim$ is an equivalence relation over $G$.

We first prove some results before giving the proof of Theorem~\ref{mainth-or}.

\begin{lemma}\label{olem1}
Let $G$ be a group such that $|G|$ is divisible by at
least two distinct primes.
Let $x$ and $y$ be two distinct elements of $G$. Then
$x\sim y$ in $\mathcal{S}(G)$ if and only if one of the
following occurs:

\noindent {\rm (i)} $o(x)=o(y)$.

\noindent {\rm (ii)} $\{o(x),o(y)\}=\{1,\exp(G)\}$.

\noindent {\rm (iii)} $\{o(x),o(y)\}=\{p^m,p^n\}$ and $p^nq\notin \pi_e(G)$, where $p,q$ are two distinct primes and $m,n$ are two positive integers with $m>n$.
\end{lemma}
\proof
By the definition of an order supergraph,
the proof of the sufficiency is straightforward.
We next prove the necessity. Suppose that $x\sim y$ in $\mathcal{S}(G)$. Assume that $o(x)\ne o(y)$.
Suppose that one of $x$ and $y$ is $e$. Without loss of generality, let $x=e$. Then $N[y]=G$. Since $|G|$ is divisible by at least two distinct primes, we have that $o(y)$ is not a prime power.
It follows from $N[y]=G$ that $\exp(G)\mid o(y)$.
Also, as $o(y)$ divides $\exp(G)$, we actually have that
$\exp(G)=o(y)$, as desired.

Suppose, in the following, that $e\notin \{x,y\}$.
We claim that if $o(x)$ is not a prime power, then $o(x)\mid o(y)$.  In fact,
let $q^t\mid o(x)$ and $q^{t+1}\nmid o(x)$, where $q$ is a prime. It follows that there exists $a\in G$ such that $o(a)=q^t$, and so $a\in N[y]$. Note that $o(x)$ is not a prime power. Let $r\ne q$ be a prime divisor of $o(x)$.
It follows that there exists an element of order $r$ such that it belongs to $N[x]=N[y]$, which implies that $o(y)$ is not a power of $q$. As a result, we have $a\ne y$. It follows that $q^t\mid o(y)$, and so $o(x)\mid o(y)$. Thus, the claim is valid.
We conclude that if $o(x)$ is not a prime power, then $o(y)$ is also not a prime power, it follows from the above claim that $o(y)=o(x)$, a contradiction.
So, we may assume that $o(x)=p^m$ and $o(y)=p^n$ for some prime $p$ and two distinct positive integers $m,n$. Without loss of generality, we may assume that $m>n$. Suppose, to the contrary, that there exists an element $z$ in $G$ such that
$o(z)=p^nq$ for some prime $q\ne p$. Then $z\in N[y]$, and so $z\in N[x]$. It follows that $p^m\mid p^nq$, contrary to $m>n$. Thus, the necessity follows.
$\qed$

The following result is immediate by Lemma~\ref{olem1}.

\begin{cor}\label{ocor1}
Let $x,y\in G$ with $\{o(x),o(y)\}=\{p^m,p^n\}$, where $p$ is a prime and $m,n$ are positive integers with $m>n$.
Then $x\sim y$ if and only if
$p^nq\notin \pi_e(G)$ for any prime $q\ne p$.
\end{cor}

For some elements $a_1,a_2,\ldots,a_k$ of $G$, if
$o(a_1)\mid o(a_2)\mid\cdots\mid o(a_k)$ and $o(a_i)\ne o(a_j)$ for any two indices $1\le i < j \le k$, then $\{a_1,a_2,\ldots,a_k\}$ is called a {\em proper order chain} of $G$.

\begin{lemma}\label{olem2}
If $C$ is a clique of $\mathcal{R}_{\mathcal{S}(G)}$, then $C$ is a proper order chain of $G$.
\end{lemma}
\proof
Notice that $o(x)\ne o(y)$ for each two distinct $x,y\in C$.
We proceed by induction on the size of $C$. If $|C|=2$, the desired result follows.
Assume inductively that the result holds
for cliques of size $n$.
Let $C=\{a_1,a_2,\ldots,a_n,a_{n+1}\}$. Then, without loss of generality, we may assume that
$o(a_1)\mid o(a_2)\mid\cdots\mid o(a_n)$ and $\{a_1,a_2,\ldots,a_n\}$ is a proper order chain.
If $o(a_{n+1})\mid o(a_1)$, then the desired result follows.
As a result, we may assume that $o(a_1) \mid o(a_{n+1})$.
Let
$$k=\max\{i:o(a_i)\mid o(a_{n+1})\}.$$
If $k=n$, then the required result follows. Otherwise, we must have $o(a_k)\mid o(a_{n+1})\mid o(a_{k+1})$,
as desired.
$\qed$

A graph is called a {\em tree} if it is connected and has no cycles. A graph is called a {\em star} if it is a tree on $n$  vertices with one vertex having degree $n-1$ and the other
$n-1$ vertices having degree $1$.

\begin{thm}\label{mth-or}
Let $G$ be a group of order $n$. Then
$$\omega(\mathcal{R}_{\mathcal{S}(G)})
=\left\{
                                  \begin{array}{ll}
                                  1, & \hbox{if $G$ is a $p$-group;} \\
                                  \Omega(n), & \hbox{if $G$ is a cyclic group with at least two distinct prime divisor;} \\
                                    2, & \hbox{if $G$ is a CP-group with at least two distinct prime divisors;} \\
                                    \lambda_G+1, & \hbox{otherwise,}
                                  \end{array}
                                \right.
$$
where $\lambda_G=\max\{\Omega(m): m\in \pi_e(G) \text{ and $m$ is not a prime power}\}$.
\end{thm}
\proof
Note that $\mathcal{S}(G)$ is complete if and only if $G$ is a $p$-group. Thus, if $G$ is a $p$-group, then $\mathcal{R}_{\mathcal{S}(G)}$ has order $1$, and so $\omega(\mathcal{R}_{\mathcal{S}(G)})=1$, as desired.
Suppose now that $G$ is a cyclic group with at least two distinct prime divisors. Then it follows from \cite[Theorem 2.2]{HA} that $\mathcal{S}(G)=\mathcal{P}(G)$. Thus,
in view of \cite[Theorem 3.1]{Mas}, we have $\omega(\mathcal{R}_{\mathcal{S}(G)})=\Omega(n)$, as desired.

Suppose next that $G$ is a CP-group with at least two distinct prime divisors. Then $G$ is non-cyclic.
By Lemma~\ref{olem1}, for distinct $x,y\in G$,
we have that $x\sim y$ if and only if $o(x)=p^m$ and $o(y)=p^n$ where $p$ is a prime. It follows that
$\mathcal{R}_{\mathcal{S}(G)}$ is a star, which implies that
$\omega(\mathcal{R}_{\mathcal{S}(G)})=2$, as desired.

Finally, suppose that $G$ is a non-cyclic group with at least two distinct prime divisors and is not a CP-group.
Let $C=\{a_1,a_2,\ldots,a_t\}$ be a clique of $\mathcal{R}_{\mathcal{S}(G)}$ with
$|C|=\omega(\mathcal{R}_{\mathcal{S}(G)})$.
Then from Lemma~\ref{olem2}, it follows that $C$ is a proper order chain of $G$. Thus, without loss of generality, we may assume that $o(a_1)\mid o(a_2)\mid\cdots\mid o(a_t)$.
Note that
$$\lambda_G=\max\{\Omega(m): m\in \pi_e(G) \text{ and $m$ is not a prime power}\}.$$
In the following, we first prove
\begin{equation}\label{meq1}
|C|\le \lambda_G+1.
\end{equation}
If $o(a_t)$ is not a prime power, then it is easy to see that $|C|\le \lambda_G+1$, as desired. Now suppose that $o(a_t)=p^k$ for some prime $p$ and positive integer $k$.
If $a_{t-1}=e$, then $|C|=2< \lambda_G+1$ since $G$ is not a CP-group, as desired. As a result, we may assume that $o(a_{t-1})=p^l$ for some $1\le l <k$. Note that $N[a_t]\ne N[a_{t-1}]$. By Corollary~\ref{ocor1}, there exists $x\in G$ such that $o(x)=p^lq$ for some prime $q\ne p$. Therefore, $\{a_1,a_2,\ldots,a_{t-1},x\}$ is also a clique of $\mathcal{R}_{\mathcal{S}(G)}$, which implies that
$|C|\le \Omega(o(x))+1\le \lambda_G+1$, as desired.

On the other hand, let
\begin{equation*}\label{}
   m=p_1^{r_1}p_2^{r_2}\cdots p_h^{r_h}\in \pi_e(G),
\end{equation*}
where $h\ge 2$, $p_1,p_2,\ldots,p_h$ are pairwise distinct primes and $r_i\ge 1$ for any $1\le i \le h$.
Take $y\in G$ with $o(y)=m$.
Now let
$T=\{e,y_1,y_2,\cdots,y_{\Omega(m)}\}$ be a subset of $\langle y\rangle$ such that
$$
\begin{array}{l}
|y_1|=p_1,|y_2|=p_1p_2,|y_3|=p_1^2p_2,|y_4|=p_1^3p_2,
\ldots,|y_{r_1+1}|=p_1^{r_1}p_2,\\
|y_{r_1+2}|=p_1^{r_1}p_2^2,|y_{r_1+3}|=p_1^{r_1}p_2^3,
\ldots,|y_{r_1+r_2}|=p_1^{r_1}p_2^{r_2},\\
|y_{r_1+r_2+1}|=p_1^{r_1}p_2^{r_2}p_3,|y_{r_1+r_2+2}|
=p_1^{r_1}p_2^{r_2}p_3^2,\ldots,|y_{r_1+r_2+r_3}|=
p_1^{r_1}p_2^{r_2}p_3^{r_3}, \\
\ldots \ldots
\\
|y_{r_1+r_2+\ldots+r_{h-1}+1}|=p_1^{r_1}p_2^{r_2}\cdots p_{h-1}^{r_{h-1}}p_h,|y_{r_1+r_2+\ldots+r_{h-1}+2}|=p_1^{r_1}p_2^{r_2}\cdots p_{h-1}^{r_{h-1}}p_h^2,\ldots,\\
|y_{r_1+r_2+\ldots+r_{h-1}+r_h-1}|=p_1^{r_1}p_2^{r_2}\cdots p_{h-1}^{r_{h-1}}p_h^{r_h-1},|y_{\Omega(m)}|=m.
\end{array}
$$
Note that $G$ is neither a $p$-group nor a cyclic group.
By Lemma~\ref{olem1}, it is easy to see that $T$ is a clique of
$\mathcal{R}_{\mathcal{S}(G)}$ with size $\Omega(m)+1$. It follows that $\mathcal{R}_{\mathcal{S}(G)}$ has a clique of size
$\lambda_G+1$. Now (\ref{meq1}) implies that $\omega(\mathcal{R}_{\mathcal{S}(G)})=\lambda_G+1$, as required.
$\qed$

\medskip

Theorem~\ref{mainth-or} follows from Theorems~\ref{Newthm0}
and \ref{mth-or}.

\section{Enhanced power graphs}

Panda {\em et al.} \cite{PDK} computed the strong metric dimensions of the enhanced power graphs of some groups, such as, dihedral groups and semi-dihedral groups.
In this section, we characterize
the strong metric dimension of the
enhanced power graph of a group (see Theorem~\ref{Eth1}).

Let $G$ be a group.
For any $g\in G$, define
$$
[g]:=\{x\in G: \langle x\rangle=\langle g\rangle\},
$$
$$
\mathcal{M}_g:=\{M\in \mathcal{M}_G: g\in M\},
$$
and
\begin{equation}\label{nne1}
\mathcal{C}(g):=
\bigcap_{M\in \mathcal{M}_g}M\setminus
\bigcup_{M\in \mathcal{M}_G \setminus\mathcal{M}_g}M.
\end{equation}
Note that $g\in \mathcal{C}(g)$ and that $\mathcal{C}(e)=\bigcap_{M\in \mathcal{M}_G}M$, because $\mathcal{M}_e=\mathcal{M}_G$.
For $x,y \in G$, denote by $\equiv$ the equivalence relation defined by  $N[x]=N[y]$ in $\mathcal{P}_E(G)$.
As stated in Section~\ref{pre}, $\equiv$ is an equivalence relation over $G$.
The $\equiv$-class containing the element $x\in G$ is denoted by $\overline{x}$.
Let $\overline{G}=\{\overline{x}:x\in G\}$.

Recall that $\mathcal{P}_E(G)$ is complete
if and only if $G$ is cyclic (see \cite[Theorem 2.4]{Bh18}).
Thus, if $G$ is a cyclic group, then $\overline{g}=\mathcal{C}(g)=G$ for any $g\in G$,
since $\mathcal{M}_G=\{G\}$ if and only if $G$ is cyclic.
Now in view of \cite[Proposition 2.3]{MS}, we have the following result, which characterizes every $\equiv$-class.

\begin{lemma}\label{2-lem1}
For every $g\in G$, we have
$\overline{g}=\mathcal{C}(g)$. In particular,
$[g]\subseteq \overline{g}$.
\end{lemma}

\begin{lemma}\label{2-lem2}
A maximal clique of $\mathcal{R}_{\mathcal{P}_{E}(G)}$ is a subset of some maximal cyclic subgroup of $G$.
\end{lemma}
\proof
By the definition of $\mathcal{R}_{\mathcal{P}_{E}(G)}$, it is easy to see that A maximal clique in $\mathcal{R}_{\mathcal{P}_{E}(G)}$ is also a clique  in $\mathcal{P}_{E}(G)$. Now \cite[Lemma 33]{acam} implies that a maximal clique in the enhanced power graph is a cyclic subgroup, so a maximal clique of $\mathcal{R}_{\mathcal{P}_{E}(G)}$ is a subset of some maximal cyclic subgroup of $G$.
$\qed$

\begin{lemma}\label{2-lem3}
If $\{x_1,x_2,\ldots,x_t\}$ is a maximal clique of $\mathcal{R}_{\mathcal{P}_{E}(G)}$, then $\bigcup_{i=1}^t\overline{x_i}$ is a maximal cyclic
subgroup of $G$.
\end{lemma}
\proof
By Lemma~\ref{2-lem2}, there exists
$\langle x\rangle\in \mathcal{M}_G$ such that
$\{x_1,x_2,\ldots,x_t\}\subseteq \langle x\rangle$.
Also, note that for any $1\le i \le t$, we have
$\langle x\rangle\in\mathcal{M}_{x_i}$.
It follows from Lemma~\ref{2-lem1} and (\ref{nne1}) that
$\overline{x_i}\subseteq \langle x\rangle$, and so
$\bigcup_{i=1}^t\overline{x_i}\subseteq \langle x\rangle$. It suffices to prove that $ \langle x\rangle\subseteq\bigcup_{i=1}^t\overline{x_i}$. Suppose, to the contrary, that there exists $y\in  \langle x\rangle$ such that
$y\notin \bigcup_{i=1}^t\overline{x_i}$. Then, similarly, we can deduce that  $\overline{y}\subseteq \langle x\rangle$.
Note that $y$ is adjacent to $x_i$ in $\mathcal{P}_{E}(G)$.
We then have that $\{x_1,x_2,\ldots,x_t,y\}$ is a clique of $\mathcal{R}_{\mathcal{P}_{E}(G)}$, this contradicts our hypothesis that $\{x_1,x_2,\ldots,x_t\}$ is a maximal clique of $\mathcal{R}_{\mathcal{P}_{E}(G)}$.
$\qed$

\begin{lemma}\label{2-lem4}
Let $x,y\in G$. Then

\noindent {\rm (i)} $N_{\mathcal{P}_{E}(G)}[x]=\bigcup_{M\in \mathcal{M}_x}M$.

\noindent  {\rm (ii)}
$x\equiv y$ if and only if $\mathcal{M}_x=\mathcal{M}_y$.
\end{lemma}
\proof (i) Taking $w\in N_{\mathcal{P}_{E}(G)}[x]$, we have that
$\langle x,w\rangle$ is cyclic, and so there exists a maximal cyclic subgroup $M$ such that $\langle x,w\rangle \subseteq M$.
As a result, $M\in \mathcal{M}_x$, which implies that $w\in M\subseteq \bigcup_{M\in \mathcal{M}_x}M$. So, $N_{\mathcal{P}_{E}(G)}[x]\subseteq\bigcup_{M\in \mathcal{M}_x}M$. On the other hand, for any $z\in \bigcup_{M\in \mathcal{M}_x}M$, we have $z\in N$ for some
$N\in \mathcal{M}_x$. It follows that $\langle x,z\rangle$ is cyclic, and hence $z\in N_{\mathcal{P}_{E}(G)}[x]$. Namely,
$\bigcup_{M\in \mathcal{M}_x}M\subseteq N_{\mathcal{P}_{E}(G)}[x]$, as desired.

(ii) If $\mathcal{M}_x=\mathcal{M}_y$, then (i) implies
$N_{\mathcal{P}_{E}(G)}[x]=N_{\mathcal{P}_{E}(G)}[y]$, and so
$x\equiv y$, as desired.
For the converse, suppose that $x\equiv y$. Let $\langle g\rangle\in \mathcal{M}_x$. Then $g\in N_{\mathcal{P}_{E}(G)}[x]$ by (i). Since $N_{\mathcal{P}_{E}(G)}[x]=N_{\mathcal{P}_{E}(G)}[y]$, we have that $\langle g,y\rangle$ is cyclic.
Now from $\langle g\rangle\in \mathcal{M}_G$, it follows that $\langle g,y\rangle=\langle g\rangle$, so $\langle g\rangle\in \mathcal{M}_y$. As a result, $\mathcal{M}_x\subseteq\mathcal{M}_y$. Similarly, we also can deduce $\mathcal{M}_y\subseteq\mathcal{M}_x$.
$\qed$

\medskip

Combining Lemmas~\ref{2-lem3}, \ref{2-lem4} and Theorem~\ref{Newthm0}, we obtain the main result of this section.

\begin{thm}\label{Eth1}
Let $G$ be a group of order $n$. Then
\begin{equation*}
\begin{split}
  \sdim(\mathcal{P}_{E}(G))&= n-\max\{|\overline{M}|: M\in \mathcal{M}_G\} \\
    &=n-\max\{|S|: S\subseteq M\in \mathcal{M}_G\text{ and
for any $x,y\in S$, $\mathcal{M}_x\ne \mathcal{M}_y$}\}.
\end{split}
\end{equation*}
\end{thm}

The following result is immediate by Theorem~\ref{Eth1}.

\begin{cor}
Let $G$ be a group of order $n$. Then

\noindent {\rm (i)} $\sdim(\mathcal{P}_{E}(G))=n-1$ if and only if
$G$ is cyclic.

\noindent {\rm (ii)} If $G$ is a non-cyclic $\mathcal{P}$-group, then $\sdim(\mathcal{P}_{E}(G))=n-2$.
\end{cor}

By Theorem~\ref{Eth1}, (\ref{q4n-0}) and (\ref{q4n-1}), we determine the strong metric dimension of the enhanced power graph of a generalized quaternion group.

\begin{cor}
Let $Q_{4n}$ be the generalized quaternion group as presented in {\rm (\ref{q4n})}. Then
$\sdim(\mathcal{P}_{E}(Q_{4n}))=4n-2$.
\end{cor}

As an application of Theorem~\ref{Eth1}, we determine the
strong metric dimension of the enhanced power graph of an abelian $p$-group.

\begin{prop}
Let $G$ be a non-cyclic abelian $p$-group with order $n$ and exponent $p^m$. Then $\sdim(\mathcal{P}_{E}(G))=n-m-1$.
\end{prop}
\proof
Note that $G$ is non-cyclic. We may assume that $G=A\times B$ where $A$ is an abelian $p$-group and $B=\langle b\rangle$  with $o(b)=p^m$. Then $\langle (e,b)\rangle\cong B$ is a maximal cyclic subgroup of order $p^m$. Clearly,
\begin{equation}\label{neqq1}
\langle (e,b)\rangle=[(e,b^{p^{m}})]\cup [(e,b^{p^{m-1}})]
\cup [(e,b^{p^{m-2}})]\cup \dots \cup [(e,b^{p})]\cup [(e,b^{p^0})].
\end{equation}
Let $a\in A$ with order $p$. In the following, we prove that
for any two $0\le i<j\le m$,
\begin{equation}\label{neqq2}
\overline{(e,b^{p^{i}})}\ne \overline{(e,b^{p^{j}})}.
\end{equation}
Note that $i\le j-1 \le m-1$.
Now $o((a,b^{p^{j-1}}))=p^{m-j+1}$ and
$(e,b^{p^{j}})\in \langle(a,b^{p^{j-1}})\rangle$.
Let $M\in \mathcal{M}_G$ with $\langle(a,b^{p^{j-1}})\rangle\subseteq M$. Then
$M\in \mathcal{M}_{(e,b^{p^{j}})}$. Assume, to the contrary,
that $(e,b^{p^{i}})\in M$. Note that $o((e,b^{p^{i}}))=p^{m-i}$ and $M$ is a cyclic $p$-group.
If $m-i> m-j+1$, then
$\langle(a,b^{p^{j-1}})\rangle\subseteq
\langle(e,b^{p^{i}})\rangle$, a contradiction.
Since $0\le i<j\le m$, it follows that $m-i=m-j+1$.
This means that the order of $\langle(a,b^{p^{j-1}})\rangle$ is equal to the order of $\langle(e,b^{p^{i}})\rangle$.
Since  $(e,b^{p^{i}})\in M$ and $(a,b^{p^{j-1}})\in M$,
we obtain a contradiction as
$\langle(a,b^{p^{j-1}})\rangle\ne \langle(e,b^{p^{i}})\rangle$.

We conclude $M\notin \mathcal{M}_{(e,b^{p^{i}})}$, and so
$\mathcal{M}_{(e,b^{p^{i}})}\ne \mathcal{M}_{(e,b^{p^{j}})}$.
Now Lemma~\ref{2-lem4}(ii) implies that (\ref{neqq2}) is valid.
It follows from (\ref{neqq1}) and Lemma~\ref{2-lem1} that $\overline{\langle (e,b)\rangle}=m+1$. Also, note that the fact that a maximal cyclic subgroup of order $p^t$ has at most
$t+1$ $\equiv$-classes. Since $G$ has exponent $p^m$, we have
$\sdim(\mathcal{P}_{E}(G))=n-m-1$ by
Theorem~\ref{Eth1}.
$\qed$

\section{Reduced power graphs}

In this section, we characterize the strong metric dimension of the reduced power graph of a group.
Our main result is the following theorem.

\begin{thm}\label{main-New1}
Let $G$ be a group of order $n$. Then
$$
\sdim(\mathcal{P}_{R}(G))=\left\{
                                  \begin{array}{ll}
                                    2^k-k, & \hbox{if $G\cong \mathbb{Z}_{2^k}$, where $k\ge 1$;} \\
                                    2^{t+2}-t-1, & \hbox{if $G\cong Q_{4\cdot 2^t}$, where $t\ge 1$;}\\
                                    n-\max\{\Omega(m): m\in \pi_e(G)\}-1, & \hbox{otherwise.}
                                  \end{array}
                                \right.
$$
\end{thm}

In the following, we prove some results before giving the proof of Theorem~\ref{main-New1}.

\begin{lemma}\label{lem2}
Let $x$ and $y$ be two distinct elements of $G$. Then
$N[x]=N[y]$ in $\mathcal{P}_{R}(G)$ if and only if $G$ is isomorphic to either $\mathbb{Z}_{2^m}$ or $Q_{4\cdot 2^m}$ where $m$ is a positive integer, and $\{x,y\}=\{e,a\}$ where $a$ is the unique involution of $G$.
\end{lemma}
\proof If $G\cong \mathbb{Z}_{2^m}$, clearly, $N[e]=N[a]=G$ where $a$ is the unique involution of $G$, as desired. If $G\cong Q_{4\cdot 2^m}$, it follows from (\ref{q4n-1}) that $N[e]=N[a]=G$, where $a$ is the unique involution of $G$, as desired. Thus, the sufficiency follows.

We next prove the necessity.
Let $x$ and $y$ be distinct elements of $G$ and assume that
$N[x]=N[y]$ in the graph $\mathcal{P}_{R}(G)$.
Since $y^{-1}\in N[x]=N[y]$, it follows that $y=y^{-1}$.
Similarly $x=x^{-1}$. As $x$ and $y$ are adjacent in $\mathcal{P}_{R}(G)$, we must have that
$\{x,y\}=\{e,a\}$, where $a$ is an involution.
Observe that $N[a]=N[e]=G$. From this observation,
we deduce that $G$ must be a $2$-group and that $a$ must be the unique involution of $G$.
Now in view of Theorem~\ref{p-unique}, we have that $G$ is isomorphic to either $\mathbb{Z}_{2^m}$ or $Q_{4\cdot 2^m}$,
as wanted.
$\qed$

\begin{lemma}\label{Refer-lem}
If $\mathcal{C}$ is a clique in $\mathcal{P}_{R}(G)$, then
$\langle \mathcal{C}\rangle$ is cyclic.
\end{lemma}
\proof
We shall use induction on $|\mathcal{C}|$.
The result is trivial for $|\mathcal{C}|=2$ and so
assume that $|\mathcal{C}|>2$. Fix $x\in \mathcal{C}$.
If $\langle y\rangle\subset \langle x\rangle$ for every
$y\in \mathcal{C}\setminus\{x\}$, then $\langle \mathcal{C}\rangle\subseteq \langle x\rangle$ and so
$\langle \mathcal{C}\rangle$ is cyclic. If
$\langle x\rangle\subset \langle y\rangle$ for some
$y\in \mathcal{C}\setminus\{x\}$, then
$\langle \mathcal{C}\rangle \subseteq \langle\mathcal{C}\setminus\{x\}\rangle$.
The subgroup $\langle\mathcal{C}\setminus\{x\}\rangle$ is cyclic by our induction hypothesis, and so it follows that
$\langle \mathcal{C}\rangle$ is cyclic in this case too.
The induction argument goes through.
$\qed$

\medskip

The following result determines the  clique number  of a reduced power graph, which also was proved in \cite{RP2} by an alternative method.

\begin{lemma}\label{n-cliquen}
Let $G$ be a group. Then
$\omega(\mathcal{P}_{R}(G))=\max\{\Omega(m): m\in \pi_e(G)\}+1$.
\end{lemma}
\proof
Let $k=\max\{\Omega(m): m\in \pi_e(G)\}+1$ and let
$\{x_1,x_2,\ldots,x_t\}$ be a clique of $\mathcal{P}_{R}(G)$ with
size $\omega(\mathcal{P}_{R}(G))$. It suffices to prove $t=k$.
By Lemma~\ref{Refer-lem}, we have that  $\{x_1,x_2,\ldots,x_t\}\subseteq \langle x\rangle$ for some $x\in G$.
Now let $o(x)=m$. Note that for each two $1\le i<j \le t$,
$o(x_i)\ne o(x_j)$, and $o(x_i)\mid o(x_j)$ or $o(x_j)\mid  o(x_i)$. Also, $\{x_1,x_2,\ldots,x_t\}$ must be a clique of $\mathcal{P}_{R}(\langle x\rangle)$ with size $\omega(\mathcal{P}_{R}(\langle x\rangle))$. We deduce that
$t=\Omega(m)+1$, and so $t\le k$.

On the other hand, let $n\in \pi_e(G)$ with $k=\Omega(n)+1$ and let
\begin{equation*}\label{n}
   n=p_1^{r_1}p_2^{r_2}\cdots p_m^{r_m},
\end{equation*}
where $p_1,p_2,\ldots,p_m$ are pairwise distinct primes and $r_i\ge 1$ for any $1\le i \le m$.
Take $a\in G$ with $o(a)=n$.
Let
$T=\{e,a_1,a_2,\cdots,a_{\Omega(n)}\}$ be a subset of $\langle a\rangle$ such that
$$
\begin{array}{l}
|a_1|=p_m,|a_2|=p_m^2,\ldots,|a_{r_m}|=p_m^{r_m},\\
|a_{r_m+1}|=p_{m-1}p_m^{r_m},\ldots,|a_{r_m+r_{m-1}}|=
p_{m-1}^{r_{m-1}}p_m^{r_m}, \\
  |a_{r_m+r_{m-1}+1}|=
p_{m-2}p_{m-1}^{r_{m-1}}p_m^{r_m},
\ldots,|a_{\Omega(n)-1}|=p_1^{r_1-1}p_2^{r_2}\cdots p_m^{r_m}, |a_{\Omega(n)}|=p_1^{r_1}p_2^{r_2}\cdots p_m^{r_m}.
\end{array}
$$
Now it is easy to see that $T$ is a clique in
$\mathcal{P}_{R}(G)$ with size $\Omega(n)+1$, and so
$k\le t$.
$\qed$

\begin{lemma}\label{New-rcn}
Let $G$ be a group. Then
$$\omega(\mathcal{R}_{\mathcal{P}_{R}(G)})
=\left\{
                                  \begin{array}{ll}
                                    k, & \hbox{if $G\cong \mathbb{Z}_{2^k}$, where $k\ge 1$;} \\
                                    t+1, & \hbox{if $G\cong Q_{4\cdot 2^t}$, where $t\ge 1$;}\\
                                    \max\{\Omega(m): m\in \pi_e(G)\}+1, & \hbox{otherwise.}
                                  \end{array}
                                \right.
$$
\end{lemma}
\proof
Suppose that $G\cong \mathbb{Z}_{2^k}$ or $Q_{4\cdot 2^t}$, where $k,t\ge 1$.  Lemma~\ref{lem2} implies that $\mathcal{R}_{\mathcal{P}_{R}(G)}$ is isomorphic to the subgraph of $\mathcal{P}_{R}(G)$ obtained by deleting the vertex $e$ from $\mathcal{P}_{R}(G)$. Note that $e$ is adjacent to every non-identity element of $G$ in $\mathcal{P}_{R}(G)$. As a result, we have that $\omega(\mathcal{R}_{\mathcal{P}_{R}(G)})=\max\{\Omega(m): m\in \pi_e(G)\}$. If $G\cong \mathbb{Z}_{2^k}$, then $\max\{\Omega(m): m\in \pi_e(G)\}=\Omega(2^k)=k$, as desired. Also,
if $G\cong Q_{4\cdot 2^t}$, then by (\ref{q4n-0}), we deduce  $\max\{\Omega(m): m\in \pi_e(G)\}=\Omega(2^{t+1})=t+1$, as desired.

Suppose that $G$ is neither $\mathbb{Z}_{2^k}$ nor $Q_{4\cdot 2^t}$. By Lemma~\ref{lem2}, we have that $\mathcal{R}_{\mathcal{P}_{R}(G)}$ is equal to $\mathcal{P}_{R}(G)$, and so the desired result follows from
Lemma~\ref{n-cliquen}.
$\qed$

\medskip

Remark that $\mathcal{P}_{R}(G)$ is complete if and only if $G\cong \mathbb{Z}_2$. Thus, if $G\ncong \mathbb{Z}_2$, then $\mathcal{P}_{R}(G)$ has diameter two. Note that the strong  metric dimension of a complete graph of order $n$ is $n-1$. Thus,
combining Theorem~\ref{Newthm0} and Lemma~\ref{New-rcn},
we complete the proof of Theorem~\ref{main-New1}.

\medskip

By Theorem~\ref{main-New1} and (\ref{q4n-0}), we determine the strong metric dimension of the reduced power graph of a generalized quaternion group.

\begin{cor}
Let $Q_{4n}$ be the generalized quaternion group as presented in {\rm (\ref{q4n})}. Then
$$\sdim(\mathcal{P}_{R}(Q_{4n}))
=\left\{
                                  \begin{array}{ll}
                                    2^{t+2}-t-1, & \hbox{if $n=2^t$ for some $t\ge 1$;} \\
                                    4n-\Omega(2n)-1, &  \hbox{otherwise.}
                                  \end{array}
                                \right.
$$
\end{cor}

Clearly, for a group $G$ of order $n$, $\sdim(\mathcal{P}_{R}(G))=n-1$ if and only if $G$ is isomorphic to the cyclic group of order $2$.
As a direct application of Theorem~\ref{main-New1},
we conclude the paper by characterizing all groups $G$
whose reduced power graphs have strong  metric dimension $n-2$.

\begin{cor}
The following are equivalent for a group $G$ of order $n$:

\noindent{\rm (a)} $\sdim(\mathcal{P}_{R}(G))=n-2$;

\noindent{\rm (b)} $\mathcal{R}_{\mathcal{P}_{R}(G)}$ is a star;

\noindent{\rm (c)} $G$ is isomorphic to $\mathbb{Z}_4$, $Q_8$ or a $\mathcal{P}$-group.
\end{cor}

\bigskip

\noindent \textbf{Acknowledgements}~~We are grateful to the anonymous referee for careful reading and helpful comments.

This research was supported by the National Natural Science Foundation of China (Grant Nos. 11801441 and 61976244), the Natural Science Basic Research Program of Shaanxi (Program No. 2020JQ-761),
and the Young Talent fund of University Association for Science and Technology in Shaanxi, China (Grant No. 20190507).

\end{CJK*}

\end{document}